\documentclass[12pt]{article}
\usepackage{amsmath,amssymb,amsthm,epsfig,color}

\title{The genus of a random chord diagram is asymptotically normal.}
\author{Sergei Chmutov\,\,\,\,and Boris Pittel\thanks{
Pittel's research supported in part by NSF Grant  DMS-1101237}
\hbox{ }\\
\small 
Ohio State University\\[-0.8ex]
\small\texttt{chmutov@math.ohio-state.edu and bgp@math.ohio-state.edu}}

\def\si{\par\smallskip\noindent}
\def\bi{\par\bigskip\noindent}
\def\pr{\text{ P\/}}
\def\ex{\text{E\/}}

\def\eps{\varepsilon}

\def\part{\partial}
\def\tag#1 {\eqno(#1)}

\def\var{\text{Var\/}}
\newcommand{\ig}{\includegraphics}
\newcommand{\rb}{\raisebox}
\newcommand\risS[6]{\rb{#1pt}[#5pt][#6pt]{\begin{picture}(#4,15)(0,0)
  \put(0,0){\ig[width=#4pt]{#2.eps}} #3
     \end{picture}}}

\begin{document}
\maketitle
{\center\small Mathematics Subject Classifications: 05A16, 05C80, 57M15, 60C05, 60F05}
{\center\small Keywords: chord diagrams, random, genus, limit, distribution}
\newtheorem{Theorem}{Theorem}[section]
\newtheorem{Lemma}{Lemma}[section]
\newtheorem{Proposition}{Proposition}[section]
\newtheorem{Corollary}{Corollary}[section]
\numberwithin{equation}{section}

\begin{abstract}
Let $G_n$ be the genus of a two-dimensional surface obtained by gluing, uniformly at random,
the sides of an $n$-gon. Recently Linial and Nowik proved, via an enumerational formula due
to Harer and Zagier,  that the expected value of $G_n$  is asymptotic to 
$(n - \ln n)/2$ for $n\to\infty$.  We prove 
a local limit theorem for the distribution of $G_n$, which implies that $G_n$ is asymptotically
Gaussian, with mean $(n-\ln n)/2$ and variance $(\ln n)/4$.
\end{abstract}

\section{Introduction and main result}

In topology, it is traditional to represent a surface by gluing the sides of a polygon \cite{Ma,ZT}. The information of which side is glued to which can be encoded by a chord diagram.
For example, the classical presentation of torus by
gluing the opposite sides of a square can be presented by a chord diagram with two intersecting chords:
$$\risS{-20}{torus}{}{70}{25}{30}\qquad\qquad
\risS{-30}{torus-sq}{}{70}{0}{0}\qquad\qquad
\risS{-15}{cd22}{\put(35,34){$b$}\put(35,0){$a$}
                     \put(0,35){$a$}\put(0,-2){$b$}}{40}{0}{0}\ =
``abab" 
$$
Combinatorially, a chord diagram with $n$ chords is the same as a word of length $2n$ where each letter occurs precisely twice; think about the letters written at the end of the chords and the word is to be read off the circle counterclockwise starting with the point $(0,1)$.

Given a chord diagram $D$, the surfaces represented by it can be recovered in the following way. First, attaching a lower semi-sphere to the circle of $D$ and thickening
the chords of $D$ to narrow bands we obtain a surface with boundary:
$$D=\risS{-15}{cd34}{}{40}{0}{0}\quad \risS{-1}{totor}{}{30}{0}{0}\quad
\risS{-27}{basket-3}{}{60}{40}{25}\ .
$$
Then attaching a disc, {\it face}, to each boundary component of that surface we will get the desired closed surface. We consider only orientable surfaces here. By {\em genus} of a chord diagram we understood the genus of this surface.

In this paper we are interested in the distribution of genus $G_n$ of the chord diagram chosen uniformly at random from among all $(2n-1)!!$ such diagrams with $n$ chords. For example, for $n=3$ the genus $G_3$ takes only two values 0 and 1 according to the picture:
$$\begin{array}{c||@{\quad}c@{\quad}c@{\quad}c@{\quad}c@{\quad}c}
\mbox{genus 0}&\risS{-15}{cd35}{}{40}{20}{20}&
\risS{-15}{cd35-1}{}{40}{0}{0}&\risS{-15}{cd33}{}{40}{0}{0}&
\risS{-15}{cd33-1}{}{40}{0}{0}&\risS{-15}{cd33-2}{}{40}{0}{0} \\ \hline
\mbox{genus 1}&\risS{-15}{cd31}{}{40}{30}{20}&
\risS{-15}{cd31-1}{}{40}{0}{0}&\risS{-15}{cd31-2}{}{40}{0}{0}&
\risS{-15}{cd31-3}{}{40}{0}{0}&\risS{-15}{cd31-4}{}{40}{0}{0} \\ 
 &\risS{-15}{cd31-5}{}{40}{30}{20}&
\risS{-15}{cd34}{}{40}{0}{0}&\risS{-15}{cd34-1}{}{40}{0}{0}&
\risS{-15}{cd34-2}{}{40}{0}{0}&\risS{-15}{cd32}{}{40}{0}{0} \\ 
\end{array}
$$
Our work was inspired by a recent paper of Linial and Nowik~\cite{LN} who estimated the expected value of $G_n$:
$$
\ex[G_n]\sim \frac{n -\ln n}{2}.
$$
This estimate is implied by  a harmonic sum expression for $\ex[G_n]$ they derived from a  Harer-Zagier formula~\cite{HZ} for the bivariate generating function of  $\{c_{n,g}\}$, with $c_{n,g}$
being the number of $n$-chord diagrams of genus $g$. In principle, the Harer-Zagier formula can be used to obtain  sharp
asymptotics for higher order moments of $n+1-2G_n$. The formulas get progressively messier,
which makes a distributional analysis of $G_n$ quite problematic. From two first moments we
deduce that that the standard deviation of $G_n$ is of order $(\ln n)^{1/2}$.  Using this 
simple information  and a contour integration
formula based on the Harer-Zagier formula, we prove a local limit theorem for the distribution of $G_n$, i. e. a sharp 
asymptotic estimate for the numbers $c_{n,g}$, with $|g-\ex[G_n]| \ll (\ln n)^{7/10}$. 
As a corollary, $G_n$ is shown to be asymptotically normal, with mean $(n-\ln n)/2$ and 
standard deviation  $\sqrt{\ln n}/2$.

The Harer-Zagier formula was discovered in \cite{HZ} for the purposes of computation of the Euler characteristic of moduli spaces of complex curves. It is tightly related to matrix models of quantum gravity \cite{tH}. We recommend a remarkable book \cite{LZ} for an excellent exposition of the Harer-Zagier formula and its relation to different areas of mathematics and physics.

The topological construction of surface from a chord diagram above gives a graph with single vertex and $n$ loops embedded into a surface and dividing the surfaces into a number $F$ of cells, faces. Then the Euler characteristic of the surface is $1-n+F$. On the other hand, the Euler characteristic of an orientable surface of genus $G$ is equal to $2-2G$. Thus for the number of faces $F$ we have $F=n-2G+1$. By our result, with high probability
the genus $G_n$ of a random chord diagram is very close to  $(n-\ln n)/2$. Hence the number of faces $F_n$ typically tends to be very small, of order $\ln n$, relative to the number of edges $n$. Geometrically it means that typically  there is at least one face with a large number of sides, of order
$n/\ln n$. It would be interesting to explore the distribution of the number of sides of individual faces. For example, is there typically  just one face with that many sides, or there are several (how many) such faces? Going out on a limb, we conjecture that, analogously to cycles of a
uniformly random permutation on $[n]$,  with high probability there exist several  faces, with $\Theta(n)$ sides each.

To conclude, we should mention that the primary focus  of Linial and Nowik~\cite{LN} was a random 
{\it directed\/} $n$-chord diagram, generating a random oriented surface in a different way, for which a counterpart of the Harer-Zagier formula is unknown. In that case they used an ingenious combinatorial argument to show that
$$
\ex[G_n]=\frac{n}{2}-\Theta(\ln n).
$$

\bigskip
{\bf Acknowledgment.} We would like to thank Nathan Linial and Tahl Nowik for useful comments on the first version of the paper.

\section{Harer-Zagier formula}

In terms of the distribution $\{p_{n,g}:=\pr(G_n=g)\}$, the Harer-Zagier formula~\cite{HZ} 
is equivalent
\begin{equation}\label{1}
1+2\sum_{n,g}p_{n,g}x^{n+1}y^{n+1-2g}=\left(\frac{1+x}{1-x}\right)^y,
\end{equation}
\cite{LN}. To illustrate the power of \eqref{1}, let us compute $p_{n,n/2}=\pr(G_n=n/2)$, which is
the probability that the random surface has exactly one face, $F_n=1$. (Of course, $p_{n,n/2}=0$
for $n$ odd.) It follows from \eqref{1} that
$$
1+2\sum_np_{n,n/2}x^{n+1}=[y^1]\left(\frac{1+x}{1-x}\right)^y=\ln\frac{1+x}{1-x}.
$$
Hence, for even $n$,
$$
p_{n,n/2}=[x^{n+1}]\sum_{j\text{ odd}}\frac{x^j}{j}\Longleftrightarrow p_{n,n/2}=\frac{1}{n+1}.
$$
Observe that $1/(n+1)$ is the probability that the uniformly random permutation $\omega_{n+1}$
of
$[n+1]$ is cyclic.  More generally, for $k\ge 1$,
\begin{align}
2\pr(F_n=k)=&\,[x^{n+1}]\,\frac{1}{k!}\left[\ln (1+x)-\ln (1-x)\right]^k\notag\\
=&[x^{n+1}]\,\frac{2^k}{k!}\left(\,\sum_{j\text{ odd}}\frac{x^j}{j}\right)^k.\label{more}
\end{align}
Let ${\cal O}_{a,b}$ denote the total number of permutations of $[a]$ consisting of $b$
odd cycles. From \eqref{more} and a standard exponential identity
$$
\sum_{a,b}\frac{x^ay^b}{a!}\,{\cal O}_{a,b}=\exp\left(y\sum_{j\text{ odd}}\frac{x^j}{j}\right),
$$
it follows that
\begin{multline}\label{pF=k}
\pr(F_n=k)=2^{k-1}[x^{n+1}]\sum_a\frac{x^a}{a!}{\cal O}_{a,k}
=2^{k-1}\frac{{\cal O}_{n+1,k}}{(n+1)!}\\
=2^{k-1}\pr\bigl(\omega_{n+1}\text{ consists of }k\text{ odd cycles}\bigr).
\end{multline}
\si

We will show that the equation \eqref{1} can also be used to find sharp asymptotic expression
for the moments $\ex[G_n^k]$. However, we will use only the first two moments,  as obtaining
a limiting distribution of $G_n$ via the moments method appears to be quite problematic.
Besides we set our sights higher, on a {\it local\/} limit theorem, for which the moments-based techniques
are too crude in principle. Our main tool
is contour integration already implicit  in the derivation of Theorem 2 in~\cite{HZ} from the equation \eqref{1}. This theorem (see also~\cite[Corrolary~3.1.8]{LZ}) states that
\begin{equation}\label{5}
c_{n,g}=(2n-1)!!\, p_{n,g}= \frac{(2n)!}{(n+1)!(n-2g)!}\,\,[t^{2g}]\,\left(\frac{t/2}{\tanh (t/2)}\right)^{n+1},
\end{equation}
where $[t^{2g}](f(t))$ denotes the coefficient at $t^{2g}$ in the power series Taylor expansion of the function $f(t)$.
The contour integration we mentioned is a simple consequence of \eqref{5}:
\begin{equation}\label{3}
p_{n,g}=\frac{2^n}{2\pi i\,(n-2g)!(n+1)}
\oint\limits_{C}\frac{1}{t^{2g+1}}\left(\frac{t/2}{\tanh(t/2)}\right)^{n+1}\,dt\ ;
\end{equation}
here $C$ is a simple closed contour surrounding $0$ such that all the non-zero roots of 
$\tanh (t/2)=0$ are outside of $C$. Later on we will choose $C$ depending on $n$ which will allow as to get the desired asymptotics of $p_{n,g}$.

\section{Asymptotics of $p_{n,g}$} 

The formula \eqref{1} is perfectly tailored for
asymptotic evaluation of the {\it factorial\/} moments $\ex\bigl[(n+1-2G_n)_k\bigr]$, $k\ge 1$,
($(m)_k:=m(m-1)\cdots (m-k+1)$). Indeed, differentiating
\eqref{1} $k$ times with respect to $y$, and setting $y=1$, we get
\begin{equation*}
2\sum_n x^{n+1}\ex\bigl[(n+1-2G_n)_k\bigr]=\frac{1+x}{1-x}\left(\ln\frac{1+x}{1-x}\right)^k.
\end{equation*}
So
\begin{equation*}
\ex\bigl[(n+1-2G_n)_k\bigr]=[x^{n+1}]\,\frac{1+x}{2(1-x)}\left(\ln\frac{1+x}{1-x}\right)^k.
\end{equation*}
Using an asymptotic formula for $[x^m] (1-x)^{-\alpha}\bigl(x^{-1}\ln (1-x)^{-1}\bigr)^{\beta}$,
(Flajolet and Sedgewick~\cite{FS}, Section VI.2), 
it is straightforward to write down a series-type  asymptotic formula for 
$\ex\bigl[(n+1-2G_n)_k\bigr]$. In particular, for $k=1$,
\begin{equation}\label{ex^1}
\ex[n+1-2G_n]=\ln n +\ln 2 -\Gamma^\prime(1)+O(\ln^{-1}n).
\end{equation}
Equivalently
\begin{equation}\label{exG=}
\ex[G_n]=\frac{n}{2}-\frac{\ln n}{2}+\frac{1}{2}(1-\ln 2+\Gamma^\prime(1))+O(\ln^{-1}n).
\end{equation}
This  sharp estimate can also be obtained from the harmonic sum--type formula already obtained in
\cite{LN}. Analogously, for all $k\ge 1$
\begin{equation}\label{ex^k}
\ex\bigl[(n+1-2G_n)_k\bigr]=\ln^k n +O(\ln^{k-1} n).
\end{equation}
The relation \eqref{ex^k} implies that 
$$
\frac{n+1-2G_n}{\ln n}\to 1,
$$
{\it in probability\/}. Moreover, using \eqref{ex^1} and \eqref{ex^k} for $k=2$, we see that
$$
\var[n+1-2G_n]=O(\ln n);
$$
So, by Chebyshev's  inequality,
\begin{equation}\label{CM}
\pr\bigl\{|G_n-\ex[G_n]|\ge (\ln n)^{1/2+\eps}\bigr\}\le_b (\ln n)^{-2\eps}.
\end{equation}
More expressively, by \eqref{exG=},
\begin{equation}\label{pr,G=}
G_n=\frac{n}{2}-\frac{\ln n}{2}+o_p\bigl((\ln n)^{1/2+\eps}\bigr),\quad\forall\eps>0;
\end{equation}
the $o_p$ notation means that the remainder scaled by $(\ln n)^{1/2+\eps}$ converges
to zero in probability. The upshot of \eqref{CM}-\eqref{pr,G=} is that from now we may, and will, focus
on the generic values $g$ of $G_n$ satisfying
\begin{equation}\label{upshot}
\left|g-\frac{n}{2}-\frac{\ln n}{2}\right|\le (\ln n)^{1/2+\eps}.
\end{equation}
At the risk of belaboring the obvious, the equation \eqref{CM} is equivalent to
\begin{equation}\label{risk}
\sum_{g\text{ meets }\eqref{upshot}}\!\!\!\!\!p_{n,g}= 1 - O\bigl((\ln n)^{-2\eps}\bigr).
\end{equation}

Armed with \eqref{risk}, we will determine an asymptotic formula for $p_{n,g}$, with
$g$ in the range \eqref{upshot}.  By \eqref{3},
\begin{equation}\label{2.1}
p_{n,g}=\frac{2^{-1}}{2\pi i (n-2g)!(n+1)}\oint\limits_{C}\frac{1}{t^{2g+1}}\left(\frac{t}{\tanh(t/2)}\right)^{n+1}\,dt.
\end{equation}
Since
$$
\tanh (t/2)=\frac{e^t-1}{e^t+1},
$$
its roots are $2\pi \nu i$, $\nu=\pm 1,\pm 2,\dots$. Since the integrand in \eqref{2.1} is odd,
we seek $C$ symmetric with respect to the origin $t=0$.  One would normally consider a circular contour of radius $t=t(n,g)$, where $t(n,g)$ is an 
absolute minimum point  of
$$
f(t,g)=\frac{1}{t^{2g}}\left(\frac{t}{\tanh(t/2)}\right)^{n+1}, \quad t\ge 0.
$$
However, for $g$ in the range \eqref{upshot},  $t(n,g)$  turns out to be asymptotic to $\ln n$; so a circle of that
radius would enclose not only $t=0$, but also  lots of imaginary zeroes $2\pi \nu i$, $\nu\neq 0$, of $\tanh(t/2)$. That is, this circle is inadmissible.
Instead we will select as $C$  a thin horizontal rectangular contour; its short  vertical sides pass through the points $t=\pm \bar t$, with $\bar t\sim \ln n$, and the long horizontal sides pass through the points  $t=\pm \pi i$,
the zeroes of $\coth(t/2)$ closest to the origin $t=0$.  
$$\risS{-5}{contour}{\put(1,30){$-\bar t$}\put(227,30){$\bar t$}
  \put(114,32){$0$}\put(122,61){$\pi i$}\put(121,14){$-\pi i$}}{250}{60}{0}
$$
Observe that we confine ourselves
to the same $\bar t=\bar t(n)$ for all $g$ satisfying \eqref{upshot}.
How to choose $\bar t$ ?
Our guiding intuition is that, for some $\bar g$ in the range \eqref{upshot}, $(\bar t,\bar g)$ is
a stationary, saddle-type, point of a logarithmically-sharp approximation of $f(t,g)/(n-2g)!$.
\si

Since $n-2g\sim \ln n$ for $g$ in \eqref{upshot}, we have
\begin{equation}\label{stir}
(n-2g)! =\sqrt{2\pi (n-2g)}\left(\frac{n-2g}{e}\right)^{n-2g}\bigl(1+O((\ln n)^{-1})\bigr).
\end{equation}
So we define $u=2g$, and introduce 
\begin{multline}\label{Htu=}
H(t,u)=(n-u)\ln\frac{e}{n-u}+(n+1)\ln\frac{t}{\tanh(t/2)}-u\ln t\\
=(n-u)\ln\frac{e}{n-u}+(n+1)\ln\coth (t/2)+(n+1-u)\ln t.
\end{multline}
A stationary point  of $H(t,u)$ is a solution of
\begin{align}
H_t(t,u)=&\, -(n+1)\frac{1}{\sinh t}+\frac{n+1-u}{t}=0,\label{H_t=0}\\
H_u(t,u)=&\, \ln(n-u)-\ln t.\label{H_u=0}
\end{align}
From \eqref{H_u=0}, $t=n-u$, and \eqref{H_t=0} becomes
\begin{equation}\label{t,root}
\frac{1+t}{t}\,\sinh t=n+1.
\end{equation}
Taking logarithms of both sides of \eqref{t,root}, we easily obtain
\begin{equation}\label{bar t=}
\bar t = \ln (2n) -\frac{1}{\ln (2n)} +O(\ln^{-2} n).
\end{equation}
The corresponding value $\bar g$  is
therefore
\begin{equation}\label{bar g=}
\bar g=\frac{\bar u}{2}=\frac{n-\bar t}{2}=\frac{n-\ln n}{2} +O(1).
\end{equation}
And $\bar g$ is well within the target range \eqref{upshot}! 

Using $\bar u=n-\bar t$ and 
\eqref{bar t=}, we compute
\begin{multline}\label{Hbartbaru=}
H(\bar t,\bar u)=\bar t+\ln\bar t+(n+1)\ln\frac{e^{\bar t}+1}{e^{\bar t}-1}\\
=\ln(2n)+\ln\ln n+(n+1)\bigl(2e^{-\bar t}+O(e^{-2\bar t})\bigr)\\
=\ln(2n\ln n) +O(\ln^{-1}n).
\end{multline}
Furthermore, for $g=u/2$ satisfying \eqref{upshot}, and an intermediate $\tilde u$,
\begin{multline}\label{Hbartu=}
H(\bar t,u)=H(\bar t,\bar u)+H_u(\bar t,\bar u)(u-\bar u)+\frac{1}{2}H_{uu}(\bar t,\tilde u)(u-\bar u)^2\\
=H(\bar t,\bar u)-\frac{1}{2(n-\bar u)}(u-\bar u)^2+O\bigl(|u-\bar u|^3(n-\bar u)^{-2}\bigr)\\
=H(\bar t,\bar u)-\frac{1}{2\ln n}(u-\bar u)^2+O\bigl((\ln n)^{-\delta}\bigr);
\end{multline}
here $\delta:=-3\eps+1/2>0$  if $\eps<1/6$, which we assume from now on.
\si

Putting together \eqref{upshot} and \eqref{stir}-\eqref{Hbartu=}, we transform \eqref{2.1}
into
\begin{equation}\label{transpng0}
\begin{aligned}
p_{n,g}=&\,\frac{(\ln n)^{1/2}}{(2\pi)^{3/2}i}\,\exp\left[-\frac{(u-\bar u)^2}{2\ln n}+O\bigl((\ln n)^{-\delta}
\bigr)\right]\\
&\times \oint\limits_C\frac{1}{t}\cdot\left(\frac{\coth(t/2)}{\coth(\bar t/2)}\right)^{n+1}
\cdot\left(\frac{t}{\bar t}\right)^{n+1-u}\,dt.
\end{aligned}
\end{equation}
Since $u=2g$ is even,  the integrand  in \eqref{transpng0}
is odd, just like the one in \eqref{2.1}. Consequently, the contour integral is twice the contour integral over $C^*=C_1\cup C_2
\cup C_3$; here $C_1=\{t=-\pi i+x,\,0\le x\le\bar t\,\}$, $C_2=\{t=\bar t+iy:\,-\pi\le y\le \pi\}$,
and $C_3=\{t=\pi i+x:\,\bar t\ge x\ge 0\}$. Using the main branch of $\ln z$, i. e. with the cut
$\{z:\text{Im }z=0,\,\text{Re }z\le 0\}$,  the integrand in \eqref{transpng0}  for $t\in C^*$ can be
written as $t^{-1}e^{h(t,u)}$, where
$$
h(t,u):=(n+1)\ln\frac{\coth(t/2)}{\coth(\bar t/2)}+(n+1-u)\ln\frac{t}{\bar t}.
$$
Let us show that, asymptotically, $u$ can be 
replaced with $\bar u$, i. e. the contour integral is almost independent of $u$.
\si

On $C_2$, since $\bar t\sim\ln n$,
\begin{align*}
h(t,u)=&\,h(t,\bar u)+(u-\bar u)\ln\frac{\bar t+iy}{\bar t}\\
=&\,h(t,\bar u)+O\bigl(\bar t\,^{-1}|u-\bar u|\bigr)\\
=&\,h(t,\bar u)+O\bigl((\ln n)^{\eps-1/2}\bigr),
\end{align*}
and
\begin{equation*}
|h(t,\bar u)|\le_b ne^{-\bar t}+(n+1-\bar u)\bar t\,^{-1}=O(1).
\end{equation*}
Therefore
\begin{equation}\label{intC2=}
\oint\limits_{C_1}\frac{1}{t}e^{h(t,u)}\,dt=\oint\limits_{C_1}\frac{1}{t}e^{h(t,\bar u)}\,dt
+O\bigl((\ln n)^{\eps-3/2}\bigr).
\end{equation}
\si

On $C_3$, since $\bar t=\ln n+O(1)$,
\begin{align}
h(t,u)=&\,(n+1)\ln\frac{e^{x+i\pi}+1}{e^{x+i\pi}-1}-(n+1)\ln\frac{e^{\bar t}+1}{e^{\bar t}-1}\notag\\
&+(n+1-u)\ln\frac{x+i\pi}{\bar t};\label{htuC3}\\
=&\,(n+1)\ln\frac{e^{x}-1}{e^{x}+1}+(n+1-u)\ln\frac{x+i\pi}{\bar t} +O(1).
\end{align}
Consequently
$$
\text{Re }h(t,u)=(n+1)\ln\frac{e^{x}-1}{e^{x}+1}+(n+1-u)\ln\frac{\sqrt{x^2+\pi^2}}{\bar t}+O(1).
$$
The first order derivative of the explicit part on the RHS is
$$
\frac{n+1}{\sinh x}+(n+1-u)\frac{x}{x^2+\pi^2}\ge \frac{n+1}{\sinh \bar t}\ge \frac{1}{2}.
$$
Hence
\begin{align*}
\text{Re }h(t,u)\le&\,(n+1-u)\ln\frac{\sqrt{\bar t^2+\pi^2}}{\bar t}-\frac{\bar t-x}{2}+O(1)\\
=&\,- \frac{\bar t-x}{2}+O((n+1-u)/\bar t)+O(1)\\
=&\,-\frac{\bar t-x}{2}+O(1),
\end{align*}
as $n+1-u=O(\ln n)$. Thus
$$
\bigl|e^{h(t,u)}\bigr|=e^{\text{Re }h(t,u)}\le_b e^{-(\bar t-x)/2},
$$
and picking $\gamma\in (0,1)$,
\begin{equation}\label{xsmall}
\left|\,\,\oint\limits_{t\in C_3:\,\bar t -\bar t^{\gamma}\ge x\ge 0}\frac{1}{t}\,e^{h(t,u)}\,dt\,
\right|\le_b e^{-\bar t\,
^{\gamma}/2},
\end{equation}
for all $u$ in question, including $\bar u$.
For $t\in C_3$ with $x=\text{Re }t\ge \bar t-\bar t\,^{\gamma}$,
\begin{align*}
\text{Re }(h(t,u)-h(t,\bar u))=&\,(\bar u-u)\ln\frac{\sqrt{x^2+\pi^2}}{\bar t}\\
=&\,O\bigl[(\ln n)^{1/2+\eps}\,\bar t\,^{\gamma-1}\bigr]=O\bigl((\ln n)^{-\sigma}\bigr),
\end{align*}
where $\sigma:=1/2 - \eps-\gamma>0$, if we choose $\gamma\in (0,1/2-\eps)$, which we do!
In that case
\begin{equation}\label{xlarge}
\begin{aligned}
\left|\,\,\oint\limits_{t\in C_3:\,\bar t \ge x\ge \bar t -\bar t\,^{\gamma}}\frac{1}{t}\,
\bigl[e^{h(t,u)}-e^{h(t,\bar u)}\bigr]\,dt\,\right|\le_b&\,(\ln n)^{-\sigma}\int_{\bar t-\bar t\,^{\gamma}}
^{\bar t }\frac{1}{x}\,e^{-(\bar t-x)/2}\,dx\\
\le_b\,(\ln n)^{-\sigma-1}
\end{aligned}
\end{equation}
Combining \eqref{xsmall} and \eqref{xlarge}, we obtain
\begin{align}
\oint\limits_{C_3}\frac{1}{t}\,e^{h(t,u)}\,dt-\oint\limits_{C_3}\frac{1}{t}\,e^{h(t,\bar u)}\,dt=&\,
O\bigl(e^{-\bar t\,^{\gamma}}+(\ln n)^{-\sigma-1}\bigr)\notag\\
=&\,O\bigl((\ln n)^{-\sigma-1}\bigr).\label{IntC3=}
\end{align}
The same argument yields
\begin{equation}\label{IntC1=}
\oint\limits_{C_1}\frac{1}{t}\,\bigl[e^{h(t,u)}-e^{h(t,\bar u)}\bigr]\,dt
=\,O\bigl((\ln n)^{-\sigma-1}\bigr)=O\bigl((\ln n)^{-3/2+\eps+\gamma}\bigr).
\end{equation}
Combining \eqref{intC2=}, \eqref{IntC3=} and \eqref{IntC1=},  we conclude that
\begin{equation}\label{IntC=}
\begin{aligned}
\oint\limits_{C}\frac{1}{t}\,e^{h(t,u)}\,dt=&\,2\int\limits_{C^*}
\frac{1}{t}\,e^{h(t,u)}\,dt
=\,I_n +O\bigl((\ln n)^{-3/2+\eps+\gamma}\bigr);\\
I_n:=&\,2\int\limits_{C^*}
\frac{1}{t}\,e^{h(t,\bar u)}\,dt.
\end{aligned}
\end{equation}
(We cannot write $I_n$ as the contour integral over the whole $C$ since $\bar u$
may not be an (even) integer.) 
\bi

The rest is short. Using \eqref{IntC=}, we rewrite \eqref{transpng0} as follows:
\begin{equation}\label{transpng*}
\begin{aligned}
p_{n,g}=&\,\frac{(\ln n)^{1/2}}{(2\pi)^{3/2}i}\,\exp\left[-\frac{(u-\bar u)^2}{2\ln n}+O\bigl((\ln n)^{-\delta}
\bigr)\right]\\
&\times\bigl[I_n+O\bigl((\ln n)^{-3/2+\eps+\gamma}\bigr)\bigr].
\end{aligned}
\end{equation}
Summing this expression over $g=u/2$ in the range \eqref{upshot} and using \eqref{risk}, we
get
\begin{multline}\label{sum*}
1-O\bigl((\ln n)^{-2\eps}\bigr)=\,\left[1+O\bigl((\ln n)^{-\delta}\bigr)\right]\,
\frac{(\ln n)^{1/2}}{(2\pi)^{3/2}i}\,I_n\\
\times \sum_{u:\,g=u/2\text{ meets }\eqref{upshot}}
\!\!\!\exp\left[-\frac{(u-\bar u)^2}{2\ln n}\right] +O(R_n),
\end{multline}
where
$$
R_n:=(\ln n)^{-1+\eps+\gamma}\sum_{u\in \mathbb Z}\exp\left[-\frac{(u-\bar u)^2}{2\ln n}\right].
$$
Recall also that
$$
\eps<1/6,\quad \delta=1/2-3\eps,\quad \eps+\gamma<1/2.
$$
Recognizing the sum in the definition of $R_n$ as a Riemann sum for $\int_{\mathbb R}e^{-x^2/2}dx$
times $\sqrt{\ln n}$, we see that
\begin{equation}\label{Rn=}
R_n\le_b (\ln n)^{-1/2+\eps+\gamma}.
\end{equation}
Likewise the sum in \eqref{sum*}, (with $u$ running through even integers),  equals
\begin{multline}\label{starsum=}
\frac{\sqrt{\ln n}}{2}\left(\int_{|z|\le 2(\ln n)^{\eps}}e^{-z^2/2}\,dz+O\bigl((\ln n)^{-1/2}\bigr)\right)\\
=\frac{\sqrt{2\pi \ln n}}{2}\left[1+O\bigl((\ln n)^{-1/2}\bigr)\right].
\end{multline}
Plugging the estimates \eqref{Rn=} and \eqref{starsum=} into \eqref{sum*} we arrive at
\begin{equation}\label{In=}
I_n=\frac{4\pi i}{\ln n}\left[1+O\bigr((\ln n)^{-\sigma}\bigr)
\right],
\end{equation}
where
$$
\sigma=\min\{2\eps, 1/2-\eps-\gamma,1/2-3\eps\}.
$$
By choosing $\gamma>0$ sufficiently small, we can get $\sigma$ arbitrarily close, from below, to
the {\it largest\/} value of $\min\{2\eps,1/2-3\eps\}=1/5$, which is attained at $\eps=1/10$. Combining  \eqref{transpng*} and \eqref{In=}, we have proved the following local limit theorem.
\begin{Theorem}\label{thm1}  Let $c_{n.g}$ denote the total number of chord diagrams of 
genus $g$. Then
\begin{equation}\label{transpng**}
\frac{c_{n,g}}{(2n-1)!!}=p_{n,g}=\frac{1+O((\ln n)^{-1/5+\alpha})}{\sqrt{2\pi (\ln n)/4}}\exp\left[-\frac{(g-\bar g)^2}{2(\ln n)/4}\right],
\end{equation}
uniformly for $g$ satisfying
$$
|g-\bar g|\le (\ln n)^{7/10-\alpha},\quad(\bar g=(n-\ln n)/2+O(1)),
$$
for $\alpha>0$, arbitrarily close to zero. As a corollary,  $G_n$ is asymptotically normal with mean $(n-\ln n)/2$  and variance $(\ln n)/4$.
\end{Theorem} 
\bi

\end{document}